\RequirePackage{fix-cm}
\documentclass[smallextended]{svjour3}
\smartqed

\usepackage{graphicx,latexsym,url}
\usepackage{amsmath}\usepackage{amssymb}\usepackage{amsfonts}\usepackage{amsthm}
\usepackage{color}
\usepackage[english]{babel}
\usepackage[latin1]{inputenc}
\usepackage{times}
\usepackage[T1]{fontenc}

\numberwithin{equation}{section}

\newcommand{\R}{{\rm I\!R}}
\newcommand{\pspan}{\mbox{p-span}}
\newcommand{\conv}{\mathrm{conv}}

\journalname{Optimization Letters}

\begin{document}

\title{On the Cardinality of Positively Linearly Independent Sets}
\author{W. Hare \and H. Song}
\institute{W. Hare \at
              Mathematics, University of British Columbia, 3333 University Way, Kelowna, BC, Canada.  \\
	      \email{warren.hare@ubc.ca}\\
    \and
            H. Song \at
              Mathematics, University of Calgary, 2500 University Drive NW, Calgary, AB, Canada. \\
              \email{hasong@ucalgary.ca}
}
\date{Received: date / Accepted: date}

\maketitle

\begin{abstract}  Positive bases, which play a key role in understanding derivative free optimization methods that use a direct search framework, are positive spanning sets that are positively linearly independent.  The cardinality of a positive basis in $\R^n$ has been established to be between $n+1$ and $2n$ (with both extremes existing).   The lower bound is immediate from being a positive spanning set, while the upper bound uses {\em both} positive spanning and positively linearly independent.  In this note, we provide details proving that a positively linearly independent set in $\R^n$ for $n \in \{1, 2\}$ has at most $2n$ elements, but a positively linearly independent set in $\R^n$ for $n\geq 3$ can have an arbitrary number of elements.
\end{abstract}

{\bf Keywords:} Positive Linear Independence, Positive Basis, Convex Analysis

\section{Introduction}

First introduced in 1954 \cite{Davis54}, the notion of a positive basis plays a key role in understanding derivative free optimization methods that use a direct search framework (a.k.a.\ pattern search methods) \cite{LewisTorczon1996,CoopePrice2002,KoldaLewisTorczon2003,Regis2015}.  In loose terms, a positive basis is a set that provides some directional information into every half-space, without including unnecessary vectors.

Pattern search methods seek a minimizer of $f$ by evaluating $f$ over an ever shrinking `pattern':
	\begin{quote}
	Given $f$ and a incumbent solution $x^k$. \\
	-- Evaluate $f(x^k + \alpha v^i)$ for each $v^i \in \mathcal{V} = \{v^1, v^2, ..., v^m\}$.\\
	-- If a lower function value is found, then update $x^k$.\\
	-- If no lower function value is found, then shrink $\alpha$.
	\end{quote}
If $f \in \mathcal{C}^1$, $\nabla f(x^k) \neq 0$, and the pattern forms a positive basis, then at least one vector in the pattern is a descent direction for $f$ at $x^k$.  Thus, under these conditions, once $\alpha$ is sufficiently small, descent must occur.  Formalization of this idea has been used to prove convergence in a number of direct search/pattern search methods \cite{LewisTorczon1996,CoopePrice2002,KoldaLewisTorczon2003,AudetDennis2006} (and many others).

Mathematically, to define positive basis we begin with the positive span: given a set of vectors $\mathcal{V}=\{v^1, v^2, ..., v^m\}\subseteq \R^n$, we define the positive span of $\mathcal{V}$ by
    $$\pspan(\mathcal{V}) = \left\{\sum_{i=1}^m\alpha_i v^i:\alpha_i \geq 0 ~\mbox{for all}~i=1, 2, ..., m\right\}.$$
We now define $\mathcal{V}$ to be a positive basis if it satisfies two properties
    \begin{enumerate}
        \item[i)] [Positive Spanning] $\pspan(\mathcal{V}) = \R^n$, and
        \item[ii)] [Positive Linear Independence] $v^i\notin \pspan(\mathcal{V}\backslash \{v^i\})$ for all $i=1, 2, ..., m$.
    \end{enumerate}
Unlike bases of $\R^n$ (which always have exactly $n$ elements), a positive basis of $\R^n$ has between $n+1$ and $2n$ elements.  The lower bound of $n+1$ comes from the following fact.

\begin{proposition}\cite[Thm 3.7]{Davis54} If $\mathcal{V} =\{v^1, v^2, ..., v^m\}$ is a positive spanning set of $\R^n$, then $\mathcal{V^\ominus} =\{v^1, v^2, ..., v^{m-1}\}$ is a spanning set of $\R^n$. \end{proposition}

Thus, to satisfy positive spanning condition, a set must contain at least $n+1$ vectors.

The upper bound on the cardinality of a positive basis is more subtle to obtain.  The original work of Davis includes a proof of this fact \cite[Thm 6.7]{Davis54}, but the derivation is rather technical.  Alternate proofs appear in the works of Shephard \cite{Shephard71} and Audet \cite{Audet2011}.

It is important to note that all three proofs use {\em both} positive spanning and positive linear independence to provide the upper bound.  As such, the property of positive linear independence is less clear in its implications on cardinality.  For example, any set consisting of exactly one vector is always positive linearly independent -- including the set $\mathcal{V}=\{0\}$.  (Of course, any positive linearly independent set with $2$ or more vectors cannot contain the $0$ vector.)  Based on a knowledge of positive bases, one might conjecture that positively linearly independent sets have at most $2n$ elements.  This is true in $\R^1$ and $\R^2$, but false beyond that.  In this note, we provide details proving that a positively linearly independent set in $\R^n$ for $n \in \{1, 2\}$ has at most $2n$ elements, but a positively linearly independent set in $\R^n$ for $n\geq 3$ can exceed this bound.  Proofs of the first two statements appear in Section \ref{lowdim}.  In Section \ref{highdim} we provide two examples.  The first is a straightforward example of a positively linearly independent set $\mathcal{V}$ with
$$|\mathcal{V}|={n \choose \lfloor {n \over 2} \rfloor }$$
This example requires no tools other than basic linear algebra.  The second example develops a positively linearly independent set in $\R^n$ with an arbitrary number of elements.  This example requires a small knowledge of convex analysis, but is still quite accessible to a general audience.

We conclude the introduction with a lemma that allows us to simplify notation in all future proofs and examples.

\begin{lemma}\label{lem:normalize} Let $\mathcal{V} = \{v^1, v^2, \ldots, v^m\} \subseteq \R^n$ with $v^i \neq 0$ for all $i$.  Then $\mathcal{V}$ is positively linearly independent if and only if $\widehat{\mathcal{V}} = \{{v^1 \over \|v^1\|}, {v^2 \over \|v^2\|}, \ldots, {v^m \over \|v^m\|}\}$ is positively linearly independent. \end{lemma}

\proof  Note that
    $$v^k=\sum_{i \ne k} \alpha_i v^i \quad \iff \quad {v^k \over \|v^k\|}=\sum_{i \ne k} \beta_i{v^i \over \|v^i\|},$$
where $\alpha_i \geq 0$ and $\beta_i = \alpha_i {\|v^i\| \over \|v^k\|} \geq 0$ for any $i \ne k$. Thus, ${\mathcal{V}}$ is positively linearly independent if and only if $\widehat{\mathcal{V}}$ is positively linearly independent. \qed

\section{Dimensions $\R^1$ and $\R^2$}\label{lowdim}

\begin{proposition} The maximal cardinality of a positively linearly independent set in $\R^1$ is 2. \end{proposition}

\proof Suppose $\mathcal{V} =\{v^1, v^2, v^3\}$. If $0 \in \mathcal{V}$, then $\mathcal{V}$ is not positively linearly independent. If $0 \notin \mathcal{V}$, then, by Lemma \ref{lem:normalize}, without loss of generality $v^i \in \{-1, 1\}$ for each $i$. By the pigeon hole principle, at least two vectors $i$ and $j$ are the same, and $\mathcal{V}$ is not positively linearly independent.  Thus, if $\mathcal{V}$ is positively linearly independent, then $|\mathcal{V}| \leq 2$. \qed

\begin{proposition} \label{thm:R2}The maximal cardinality of positively linearly independent sets in $\R^2$ is 4. \end{proposition}

\proof  For eventual contradiction, suppose $\mathcal{V} =\{v^1, v^2, \ldots, v^5\}\subseteq \R^2$ is positively linearly independent.  Without loss of generality, assume $\|v^i\|=1$ for $i=1, 2, ..., 5$.  Clearly, $v^1 \neq v^i$ for any $i \neq 1$.  Moreover, if $v^1 = -v^i$ for some $i$, then $v^k \neq -v^1$ for any $k\neq i$.  As such, there are at least two vectors, say $v^1$ and $v^2$, that are linearly independent.  Examining, $\mathcal{V}\setminus\{v^1, v^2\}$, and repeating the same arguments, we may assume that $v^3$ and $v^4$ are also linearly independent.  Therefore, $\{v^1, v^2\}$ and $\{v^3, v^4\}$ are both bases of $\R^2$.

As $\{v^1, v^2\}$ is a basis of $\R^2$, there exists $\alpha \in \R^2$ such that $v^3 = \alpha_1 v^1 + \alpha_2 v^2$.  If $\alpha_1 \geq 0$ and $\alpha_2 \geq 0$, then we have contradicted the positive linear independence of $\mathcal{V}$.  If $\alpha_1 > 0$ and $\alpha_2 \leq 0$, then $\alpha_1 v^1 = v^3 - \alpha_2 v^2$, which also contradicts positive linear independence.  Similarly,  $\alpha_2 > 0$ and $\alpha_1 \leq 0$ fails.  Hence, $\alpha \leq 0$.

Using the same arguments, and $\{v^3, v^4\}$ is a basis of $\R^2$, there exists $\beta \leq 0$ and $\gamma \leq 0$ such that $v^4 = \beta_1 v^1 + \beta v^2$ and $v^5 = \gamma_1 v^3 + \gamma_2 v^4.$ Substitution now yields
	$$v^5 = (\gamma_1\alpha_1 + \gamma_2\beta_1)v^1 + (\gamma_1\alpha_2 + \gamma_2\beta_2)v^2$$
which shows that $v^5 \in \pspan(\mathcal{V}\backslash \{v^5\})$, so $\mathcal{V}$ is not positively linearly independent.\qed

On a final note for this section, it is clear that the maximum cardinalities derived above can be achieved, as constructing an example with $|\mathcal{V}| = 2n$ is trivial.

\section{Dimension $\R^n$ for $n \geq 3$}\label{highdim}

Our first example provides a family of positively linear independent sets whose cardinality grows exponentially with dimension.

\begin{example}
In $\R^n$, let $\mathcal{V}$ be the set of vectors of all the nontrivial permutations with $\lfloor {n \over 2} \rfloor$ 1's and $n-\lfloor {n \over 2} \rfloor=\lceil {n \over 2} \rceil$ 0's,
	\[\mathcal{V} = \left\{ v \in \{0, 1\}^n ~:~ \sum_{i=1}^n v_i = \left\lfloor {n \over 2}\right\rfloor \right\},\]
then $|\mathcal{V}|={n \choose \lfloor {n \over 2} \rfloor}$, and $\mathcal{V}$ is positively linearly independent.
\end{example}
\proof Suppose $v^k \in \mathcal{V}$ and there exist nonnegative numbers $\alpha_i$ such that $v^k=\sum_{i\neq k}\alpha_i v^i.$ As $v^k \neq 0$, there must exist some $j$ with $\alpha_j > 0$.  Since $v^k \neq v^j$, and $v^k, v^j \in \mathcal{V}$, there must exist some index $l$ with $v^k_l = 0$ and $v^j_l=1$.  This yields the contridiction
	$$0 = v^k_l = \sum_{i\neq k} \alpha_i v^i_l \geq \alpha_j v^j_l = \alpha_j > 0.$$
Hence, $\mathcal{V}$ is positively linearly independent. The proof of the cardinality is trivial.  \qed

Our second example requires some basic tools and definitions from convex analysis.

\begin{definition}\label{df:convex}  A set $C \subseteq \R^n$ is convex, if given any $x, y \in C$ and any $\theta\in[0,1]$, the point $z=\theta x+(1-\theta)y\in C$.
\\
A set $C \subseteq \R^n$ is strictly convex, if given any $x, y \in C$ with $x \neq y$ and any $\theta\in(0,1)$, the point $z=\theta x+(1-\theta)y\in \mathrm{int}(C)$. \end{definition}

\begin{definition} The convex hull of a set $S$, denoted $\conv(S)$, is the smallest convex set that contains $S$. \end{definition}

\begin{example} The set $B_1 = \{ x \in \R^n ~:~ \|x\| \leq 1\}$ is strictly convex. \end{example}

\begin{example}\cite[Thm 2.27]{RockafellarWets1998}\label{ex:conv} Let $Y = \{y^1, y^2, ..., y^m\} \subseteq C$.  Then
    $$\conv(Y) = \left\{ x = \sum_{i=1}^m \alpha_i y^i ~:~ \alpha_i \geq 0, ~\sum_{i=1}^m \alpha_i = 1\right\}.$$\end{example}

Combining definition \ref{df:convex} and example \ref{ex:conv}, we have the following classic lemma.

\begin{lemma}\label{fact:interior} Suppose $C$ is strictly convex.  Let $Y = \{y^1, y^2, ... y^m\} \subseteq C$ and $\bar{x} \in \conv(Y)$.  Then
    $$\mbox{either}~x=y^i~\mbox{for some}~i \quad \mbox{or} \quad x \in \mathrm{int}(C).$$\end{lemma}

We now prove the main result.  

\begin{theorem}  \label{thm:Rn}For $n \geq 3$ a positively linearly independent set in $\R^n$ may contain an arbitrary number of vectors.\end{theorem}

\proof We shall provide an example in $\R^3$.  Extending to $\R^n$ can be trivially accomplished by adding $0$ elements in the extra dimensions.

Let $S=\{(x,y):x^2+y^2 = 1\}$. Select any distinct $m$ points in $S$,
	$$\{(x^1, y^1), (x^2, y^2), \ldots, (x^m, y^m)\} \subset S,$$
and define
    $$\mathcal{V}=\{v^1, v^2, \ldots, v^m\}=\{(x^1, y^1, 1), (x^2, y^2, 1), \ldots, (x^m,y^m,1)\}.$$
We claim that $\mathcal{V}$ is positively linear independent.

Suppose there exists $v^k\in \pspan(V\backslash \{v^k\})$. By reordering, without loss of generality, assume $k=m$, so
    $$v^m=\sum_{i=1}^{m-1} \alpha_i v^i, \quad \mbox{with}~~\alpha_i \geq 0.$$
Noting that $v^i_3 = 1$ for all $i$, we must have
    $$1 = v_3^m=\sum_{i=1}^{m-1} \alpha_i v_3^i = \sum_{i=1}^{m-1} \alpha_i.$$
This implies
    $$(x^m,y^m)=\sum_{i=1}^{m-1} \alpha_i(x^i,y^i), ~~\sum_{i=1}^{m-1} \alpha_i=1, ~~\alpha_i \geq 0.$$
I.e., $(x^m, y^m) \in \conv(\{(x^1, y^1), (x^2, y^2), \ldots, (x^{m-1}, y^{m-1})\})$.  As $(x^m, y^m) \neq (x^i, y^i)$ for any $i < m$, Fact \ref{fact:interior} implies $(x^m, y^m) \in \mathrm{int} (B_1)$, contradicting $(x^i, y^i) \in B_1 \setminus  \mathrm{int} (B_1)$ for all $i$.\qed

It is interesting to note that the example in the proof of Theorem \ref{thm:Rn} does not actually require $m$ to be finite.  Indeed, in $\R^3$ it is possible to create an uncountable set of vectors $\mathcal{V}$ such that  $v \notin \pspan(\mathcal{V}\backslash \{v\})$ for all $v \in \mathcal{V}$.  Moreover, while the example uses the unit sphere to create the set, it is clear that $S=\{(x,y):x^2+y^2 = 1\}$ could be replaced by the boundary of any compact strictly convex set.

Finally, we note that an alternate (but similar) proof of Theorem \ref{thm:Rn} was independently developed by Regis and presented in \cite[Thm 3.4]{Regis2015}.

\subsubsection*{Acknowledgements}

This research was partially funded by the Natural Sciences and Engineering Research Council (NSERC) of Canada, Discover Grant \#355571-2013, and by the Pacific Institute for the Mathematical Sciences (PIMS), ``Optimization: Theory, Algorithms, and Applications'' Collaborative Research Group .

The authors are indebted to Dr.\ C.\ Audet, for helpful feedback and providing the new (shorter) proof to Proposition \ref{thm:R2} which is presented in this paper.

\end{document}